\emergencystretch4mm
\input AHTOHFIE.STY
\def\Mu{{\rm M}}
\UDC{512.543.16+512.543.72}
\MSC{20f05,20f06,20e06,20e22}
\title
{THE STRUCTURE OF ONE-RELATOR RELATIVE PRESENTATIONS AND THEIR
CENTRES}
\author{Anton A. Klyachko}
\address\myAddress
\grants{\RFBR08-01-00573}

\abstract{%
Suppose that $G$ is a nontrivial torsion-free group and $w$ is a word in
the alphabet $G\cup\{x_1^{\pm1},\dots,x_n^{\pm1}\}$ such that the word
$w'\in F(x_1,\dots,x_n)$ obtained from $w$ by erasing all letters
belonging to $G$ is not a proper power in the free
group~$F(x_1,\dots,x_n)$. We show how to reduce the study of the relative
presentation $\^G=\pres<G,x_1,x_2,\dots,x_n | w=1>$ to the case $n=1$. It
turns out that any such ``$n$-variable" group $\^G$ can be obtained from
similar ``one-variable" groups by using an explicit construction similar to
wreath product. As an illustration, we prove that, for $n\ge2$, the
centre of $\^G$ is always trivial. For $n=1$, the centre of $\^G$
is also almost always trivial; there are several exceptions, and all of
them are known.

}

\s 0.
Introduction

Let $G$ be a group. A group given by a \emph{one-relator relative 
presentation over $G$} is
$$
\^G=\pres<G,x_1,x_2,\dots,x_n|w=1>\:=
G*F(x_1,x_2,\dots,x_n)/\nc w.
$$
Here $x_1,\dots,x_n$ are some letters (not belonging to $G$)
and $w$ is a word in the alphabet
$G\cup\{x_1^{\pm1},\dots,x_n^{\pm1}\}$ (such a word can be considered as
an element of the free product $G*F(x_1,x_2,\dots,x_n)$ of
$G$ and the free group with basis $x_1,x_2,\dots,x_n$). In other words,
the presentation of the group $\^G$ is obtained from a presentation
$G=\pres<A|R>$ of $G$ by adding several new generators and one
new relator:
$\^G=\pres<A\cup\{x_1,x_2,\dots,x_n\}|R\cup\{w\}>$.

Such groups $\^G$ are natural generalisations of one-relator groups 
and have been studied by many authors (see, e.g., 
[How87], [BoP92], 
[DuH93], [Met01], [Kl06b] and the references therein). 
To obtain any meaningful result, it is necessary to impose some 
constraints on the group $G$ and/or relation $w$. In this paper, we 
assume only two restrictions:
\item{(a)}
the group $G$ is torsion-free;
\item{(b)}
the word $w\in G*F(x_1,x_2,\dots)$ is such that the word $w'\in 
F(x_1,x_2,\dots)$ obtained from $w$ by erasing coefficients belonging 
to $G$ is not a proper power%
\fn{%
We say that an element $h$ of a group $H$ is a \emph{proper power} if there
exists an element $h'\in H$ and an integer $k\ge2$
such that $h=(h')^k$. In particular, the identity element \emph{is} a
proper power: $1=1^2$.
}
in the free group~$F(x_1,x_2,\dots)$.

\enditem
The same situation was considered in [Kl06a], [Kl06b], and 
[Kl07] (and also in [Kl93], [FeR96], [CR01], [Kl05], and [FoR05] in the 
case $n=1$).

The main result of this paper (Theorem 3) shows how to reduce studying  
the groups $\^G$ to the case $n=1$. It turns out that the ``$n$-variable" 
group $\^G$ can be obtained from similar ``one-variable" groups by using 
an explicit construction that involves  
\emph{free iterated amalgamated products} (see Section~3) and 
\emph{amalgamated semidirect products} (see Section~4).

Probably, this structural theorem (Theorem 3) may have 
many applications, one of which is considered in this paper. Namely, we 
study the centre of the group $\^G$.  

First, we consider the case $n=1$, i.e., the following situation.  
Let $G$ be a torsion-free group. Suppose that the group $\~G$ is obtained 
from $G$ by adding one generator and one {\it unimodular} relator, i.e., a 
relator with exponent sum one:  
$$ 
\~G=\gp{G,t\ |\ w=1}\:=(G*\gp t_\infty)/\!\nc w, 
\hbox{ where $w\equiv g_1t^{\epsilon_1}\dots g_qt^{\epsilon_q}$, 
\ \ 
$g_i\in G$,
\ \ 
$\epsilon_i\in\Z$,\ \ and\ \ $\sum\epsilon_i=1$}.  
$$ 
In this case, we 
say that the group $\~G$ is given by a \emph{unimodular relative 
presentation over $G$}. It is known that  
$\~G$ inherits some properties of the initial group $G$. In particular,
\- the abelianisations of these groups are isomorphic:
$G/[G,G]\iso \~G/[\~G,\~G]$;
\- $G$ embeds (naturally) into $\~G$
[Kl93] (see also [FeR96]); therefore, $\~G$ is nontrivial if $G$ is
nontrivial, $\~G$ is nonabelian if $G$ is nonabelian, etc.%
\footnote{**$^)$}{%
However, the natural mapping $G\to\~G$ is never
surjective, except in the case when $w\equiv gt$ [CR01].
}
\- $\~G$ (as well as $G$) is torsion-free [FoR05];
\- $\~G$ is nonsimple if $G$ is nonsimple [Kl05];
\- the Tits alternative holds for $\~G$ (i.e., $\~G$ either contains
a nonabelian free subgroup or is virtually solvable)
provided that it holds for $G$ [Kl07].

\noindent
In this paper, we establish yet another property of this kind:
\- the centre of $\~G$ is either trivial or isomorphic to the centre
of the initial group $G$.

\noindent
More precisely, we prove the following theorem.

\Th 1.
If a group $G$ is torsion-free and a word $w\in G*\gp{t}_\infty$
is unimodular, then the centre of the group $\~G=\gp{G,t\ |\ w=1}$ is
trivial, except in the following two cases:
\item{\rm1)}
$w\equiv gtg'$, where $g,g'\in G$ (and so $\~G\iso G$),
and the centre of $G$ is nontrivial;
\item{\rm2)}
the group $G$ is cyclic and $\~G$ is a one-relator group with nontrivial
centre.

One-relator groups with nontrivial centre have been well studied
([Mu64], [BaTa68], [Pi74]). The centre of each such group is infinite cyclic
(except in the case when the entire group is free abelian of rank 2; this
case is impossible in our situation). The simplest nonobvious example of
a unimodular presentation with nontrivial centre is the braid group on three
strands $\~G=B_3=\pres<g,t|gtg=tgt>$. The centre of this group is
generated by $(gt)^3$.

In the general case, any calculations in the group $\~G$ are
difficult, because the word problem in this group is not solved
(so far). For example, the natural method of finding the centre by
the formula
$$
\hbox{centre of }\~G=(\hbox{centraliser of }G)\cap(\hbox{centraliser of }t)
$$
does not work, because we can calculate none of these centralisers, so
to find the centre, we need some man{\oe}uvre.
Actually, the proof of Theorem 1 is not long, but it is heavily based on
some results from [Kl05].

Next, we pass to the ``multivariable" case. Actually, 
we consider even a more general situation.
In [Kl06a], we suggested a generalisation of the notion of unimodularity
to the case when the word $w$ is an element of the free product of a
group $G$ and any (i.e., not necessarily cyclic) group~$T$.
In this paper, we need an even more general definition.
We say that a word $w\equiv g_1t_1\dots g_qt_q\in G*T$ is
{\it generalised unimodular} if
\item{1)}
         $\prod t_i\ne1$, and the group $T$ is torsion-free;
\item{2)}
         the cyclic subgroup $\gp{\prod t_i}$ of $T$
         is a free factor of some normal
         subgroup $R=\gp{\prod t_i}*S$ of $T$;
\item{3)}
         the quotient group $T/R$ is a group with the strong unique-product
         property.

\enditem
Recall that a group $H$ is called a {\it UP-group, {\rm or a group with
the} unique-product property}, if the product $XY$ of any two finite
nonempty subsets $X,Y\subseteq H$ contains at least one element which
decomposes uniquely into the product of an element from $X$ and an element
from $Y$.
Some time ago, there was the conjecture that any torsion-free
group is UP (the converse is, obviously, true).  However, it turned out
that there exist counterexamples ([Pr88], [RS87]).

We say that a group $H$ has the {\it strong unique product property} if
the product $XY$ of any two finite
nonempty subsets $X,Y\subseteq H$ such that $|Y|\ge2$ contains at least
two uniquely decomposable elements $x_1y_1$ and $x_2y_2$ such that
$x_1,x_2\in X$,\ \ $y_1,y_2\in Y$, and $y_1\ne y_2$.

As far as we know, all known examples of UP-groups have the strong
UP-property. In particular, all right orderable groups, locally indicable
groups, and diffuse groups in the sense of Bowditch have the strong UP
property.

The main examples of generalised unimodular presentations are groups
of the form 
$$
\^G=\pres<G,x_1,x_2,\dots,x_n | w=1>,
$$ 
where the word $w\in G*F(x_1,x_2,\dots)$ is such that the word
$w'\in F(x_1,x_2,\dots)$ obtained from $w$ by erasing all coefficients
belonging to $G$ is not a proper power
in the free group~$F(x_1,x_2,\dots)$.

Indeed,
suppose that the word $w$ has the form
$w\equiv g_1x_{j_1}^{\epsilon_1}g_2x_{j_2}^{\epsilon_2}\dots
g_qx_{j_q}^{\epsilon_q}$ and $w'\in F(x_1,\dots, x_n)$ is obtained
from $w$ by erasing all coefficients:
$w'=x_{j_1}^{\epsilon_1}x_{j_2}^{\epsilon_2}\dots x_{j_q}^{\epsilon_q}$.
Consider the groups
$$
T=F(x_1,\dots,x_n)
\quad\hbox{and}\quad
T_1=\gp{x_1,\dots,x_n\ |\ w'=1}=T/\!\nc{w'}.
$$
By the Brodskii theorem [B84], if $w'$ is not a proper power in the free
group $F(x_1,\dots,x_n)$, then the group $T_1$ is locally indicable and,
hence, has the
strong UP property. By the Cohen--Lyndon theorem [CoLy63],
the element $w'$ is a primitive element of the free subgroup
$\nc{w'}$ of $T$. Thus, the word $w$, considered as an element
of the free product~$G*T$, is generalised unimodular.

In Section~5, we prove our main result, Theorem~3. As a
corollary of this structural theorem, in Section~6, we 
obtain the following fact generalising Theorem~1.

\Th 2.
Suppose that $G$ and $T$ are torsion-free groups and a cyclically
reduced word $w=g_1t_1\dots g_qt_q\in G*T$ is
generalised unimodular. Then
\item{\rm1)}
the natural mapping $G\to \^G=\pres<G,T | w=1>\:=(G*T)/\!\nc w$
is injective;
\item{\rm2)}
if the centre of $\^G$ is nontrivial and $G$
is noncyclic, then $q=1$ and either $t_1\in Z(T)$ and
$\gp{g_1}\cap Z(G)\ne1$
{\rm (in this case,
$\^G=G\zvezda_{g_1=t_1^{\hbox to 0pt{$\scriptscriptstyle{-1}$\hss}}}T$
is a free product of $G$ and $T$
with amalgamated cyclic subgroups)}
or the group $T$ is
cyclic {\rm (in this case, $T=\gp{t_1}$ and $\^G\iso G$)}.

\noindent
This theorem implies a multivariable analogue of Theorem 1.

\Corollary 1.
Suppose that $G$ is a
nontrivial torsion-free group and
a word $w\in G*F(x_1,x_2,\dots)$ is such that the word
$w'\in F(x_1,x_2,\dots)$ obtained from $w$ by erasing all coefficients
belonging to $G$ is not a proper power
in the free group~$F(x_1,x_2,\dots)$. Then
\item{\rm1)} {\rm[Kl06a]}
the natural mapping $G\to \^G=\pres<G,x_1,x_2,\dots,x_n | w=1>$
is injective;
\item{\rm2)}
if $n\ge2$,
then the centre of $\^G$ is trivial.

\Proof
If the group~$G$ is noncyclic, the assertion follows immediately 
from Theorem~2. If $G$ is cyclic, then $\^G$ is a one-relator
group with at least three generators; 
the triviality of the centres of such groups
is well known [Mu64].

\smallskip

These results on the centre of $\^G$ are not surprising. 
However, they easily implies the Kervaire--Laudenbach 
conjecture for torsion-free groups [Kl93], i.e., the nontriviality of each 
group of the form 
$$ 
\pres<H,t|w=1>, 
\ \hbox{where $H$ is a nontrivial 
torsion-free group and $w$ is any word in the alphabet 
$H\cup\{t^{\pm1}\}$}.  
$$ 
Indeed, if the group $\~H=\pres<H,t|w=1>$ is trivial, then the word $w$ 
must be unimodular (otherwise, $\~H$ admits an epimorphism onto a 
nontrivial cyclic group). Therefore, the group 
$\^H=\pres<H,t,x|w=1>=\~H*\gp x_\infty$ is centreless. This can be derived 
from either Theorem 1 (by setting $G=H*\gp x_\infty$) or Corollary 1(2) 
(by setting $G=H$). Clearly, the triviality of the centre of $\^H$ implies 
the nontriviality of $\~H$.  Thus, both Theorem 1 and Corollary 1(2) can 
be considered as strengthenings of the main result of [Kl93].  

Other theorems on the centres of one-relator
relative presentations can be found in [B81] and [Met01]. 
A different ``strong-noncommutativity" property of such
presentations was proven in [Kl06b]: if $G$ is a
nontrivial torsion-free group and $n\ge2$, then the group
$\pres<G,x_1,x_2,\dots,x_n | w=1>$ is always SQ-universal.

The author thanks James Howie and an anonymous referee for useful 
comments. 

\proclaim{Notation}\rm{
which we use is mainly standard. Note only that if $k\in\Z$, $x$ and $y$ 
are elements of a group, and $\phi$ is a homomorphism from this group into 
another group, then $x^y$, $x^{ky}$, $x^{-y}$, $x^\phi$, $x^{k\phi}$, and 
$x^{-\phi}$ denote $y^{-1}xy$, $y^{-1}x^ky$, $y^{-1}x^{-1}y$, $\phi(x)$, 
$\phi(x^k)$, and $\phi(x^{-1})$, respectively; the commutator $[x,y]$ is 
understood as $x^{-1}y^{-1}xy$.  If $X$ is a subset of a group, then 
$\gp{X}$, $\nc{X}$, and $C(X)$ denote the subgroup generated by $X$, the 
normal subgroup generated by $X$, and the centraliser of $X$, 
respectively. The centre of a group $G$ is denoted by $Z(G)$.  The symbol 
$|X|$ denotes the cardinality of a set $X$.  The letters $\Z$ and $\N$ 
denote the set of integers and the set of positive integers, 
respectively}.


\s 1.
Howie diagrams

In this section, we recall (following [Kl05]) some facts concerning 
diagrams introduced in [How83]. The only new result of this section is 
Lemma 3. 
 
Throughout this paper, the term ``surface" means ``closed oriented
two-dimensional surface".

A {\it map} $\Mu$ on a surface $S$ is a finite set of continuous mappings
$\{\mu_i\:D_i\to S\}$, where $D_i$ is a compact oriented two-dimensional
disk, called the $i$th {\it face}, or {\it cell}, of the map; the boundary
of each face $D_i$ is partitioned into finitely many intervals
$e_{ij}\subset \d D_i$, called the {\it pre-edges} of the map, by a nonempty
set of points $c_{ij}\in \d D_i$, called the {\it corners} of the map.
The images of the corners $\mu_i(c_{ij})$ and pre-edges $\mu_i(e_{ij})$
are called the {\it vertices} and {\it edges} of the map, respectively.
It is assumed that
\item{1)}
 the restriction of $\mu_i$ to the interior of each
 face $D_i$ is a homeomorphic embedding preserving orientation; the
 restriction of $\mu_i$ to each pre-edge is a homeomorphic embedding;
\item{2)}
 different edges do not intersect;
\item{3)}
  the images of the interiors of different faces do not intersect;
\item{4)}
 $\bigcup\mu_i(D_i)=S$.

\noindent
Sometimes, we interpret a map $\Mu$ as a continuous mapping
$\Mu\:\coprod D_i\to S$ from a discrete union of disks onto the surface.

The union of all vertices and edges of a map is a graph on the surface,
called the {\it $1$-skeleton}.

We say that a corner $c$ is a corner at a vertex $v$ if $\Mu(c)=v$.  There
is a natural cyclic order on the set of all corners at a vertex $v$; we
call two corners at $v$ {\it adjacent} if they are neighboring with
respect to this order.

By abuse of language, we say that a point or a subset of the surface is
contained in a face $D_i$ if it lies in the image of $\mu_i$. Similarly,
we say that a face $D_i$ is contained in some subset $X\subseteq S$ of
the surface $S$ if $\Mu(D_i)\subseteq X$.

Figure 1 presents a map on the sphere with 5 faces --- $A$, $B$, $C$,
$D$, and $E$, 18 corners --- $a_i$, $b_i$, $c_i$, $d_i$, and $e_i$, 6
vertices, 9 edges, and 18 pre-edges. Note that the number of corners always
equals the number of pre-edges and is twice the number of edges,
and the value
$$
e(S)\:=(the\ number\ of\ vertices)-(the\ number\ of\ edges)+
(the\ number\ of\ faces)
$$
does not depend on the choice of a map on the surface $S$ and is called
{\it the Euler characteristic} of this surface. The Euler characteristic
of the sphere (the only surface of our real interest in this paper)
is two.

\goodbreak
\bigskip
\centerline{\input 1.PIC}
\nobreak%
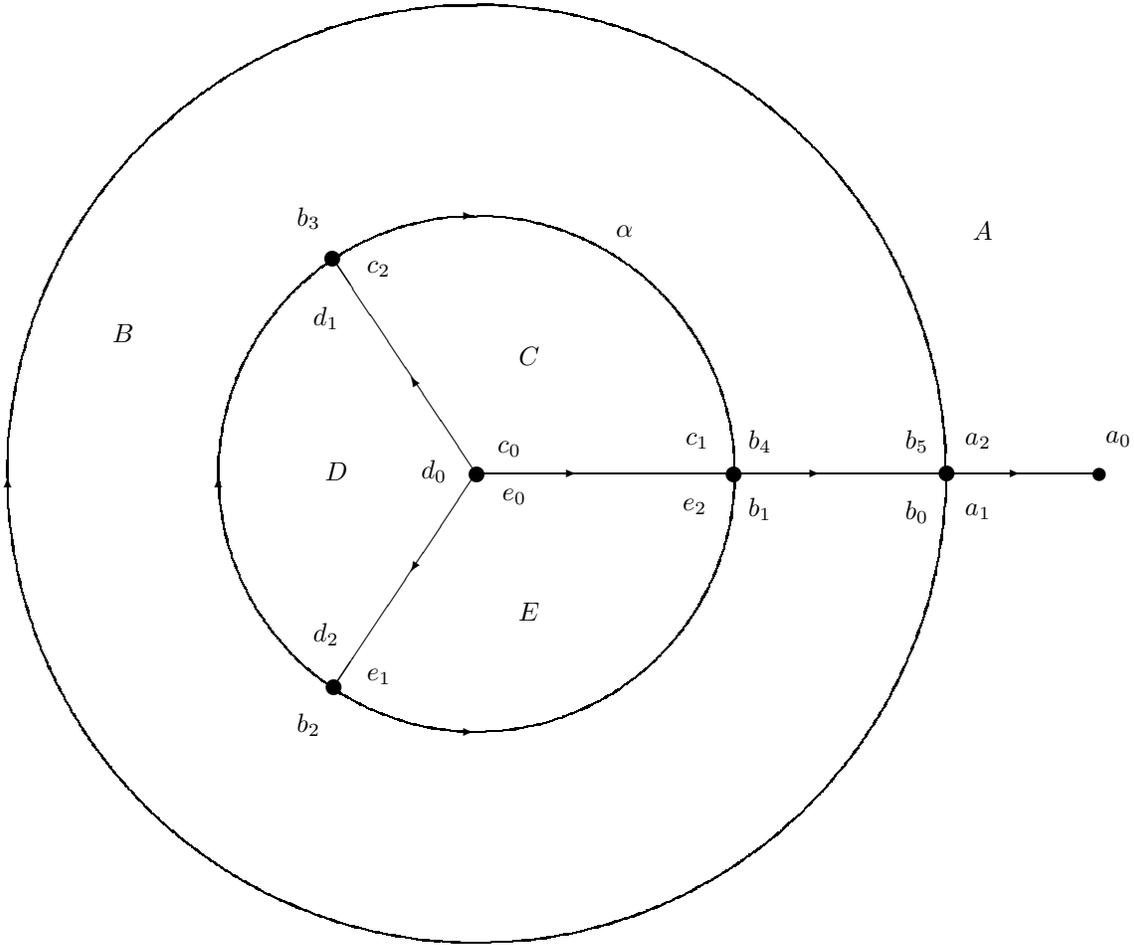
\centerline{Fig. 1}%
\goodbreak
\bigskip

Suppose that we have a map $\Mu$ on a surface $S$, the corners of the
map are labeled by elements of a group $H$, and the edges are oriented (in
the figures, we draw arrows on the edges) and labelled by elements of a
set $\{t_1,t_2,\dots\}$ disjoint from the group $H$. The label of a corner
or an edge $x$ is denoted by $\lambda(x)$.

The {\it label of a vertex} $v$ of such a map is defined by the formula
$$
\lambda(v)=\prod_{i=1}^k \lambda(c_i),
$$
where $c_1,\dots,c_k$ are all corners at $v$ listed clockwise.
The label of a vertex is an element of the group $H$ determined up to
conjugacy.

For instance, the label of the uppermost vertex in Fig. 1 
is $\lambda(b_3)\lambda(c_2)\lambda(d_1)$. 

The {\it label of a face} $D$ is defined by the formula
$$
\lambda(D)=\prod_{i=1}^k
\bigl(\lambda(\Mu(e_i))\bigr)^{\epsilon_i}\lambda(c_i),
$$
where $e_1,\dots,e_k$ and $c_1,\dots,c_k$ are all pre-edges and all
corners of $D$ listed anticlockwise, the endpoints of $e_i$ are
$c_{i-1}$ and $c_i$ (subscripts are
modulo $k$), and $\epsilon_i=\pm1$ depending on whether the homeomorphism
$e_i\mathop\to\limits^\Mu\Mu(e_i)$ preserves or reverses orientation.
Simply speaking, to obtain the label of a face, we should go around its
boundary anticlockwise, writing out the labels of all corners and edges we
meet, the label of an edge traversed against the arrow should be raised to
the power $-1$.

The label of a face is an element of the group $H*F(t_1,t_2,\dots)$ (the
free product of $H$ and the free group with basis $\{t_1,t_2,\dots\}$)
determined up to a cyclic permutation. More precisely, the right-hand side
of our formula for $\lambda(D)$ is called the {\it label of the face $D$
written starting with the pre-edge $e_1$}.

For instance, if the label of each edge in Fig. 1 is $t$, then the label 
of the face $B$ written starting with the pre-edge~$\alpha$ is 
$$ 
t\lambda(b_4)t\lambda(b_5)t^{-1}\lambda(b_0)t^{-1}\lambda(b_1) 
t^{-1}\lambda(b_2)t\lambda(b_3).
$$

Such a labelled map is called a {\it Howie diagram} (or
simply {\it diagram}) over a relative presentation
$$
K=\gp{H,t_1,t_2,\dots\ |\ w_1=1,w_2=1,\dots}
\eqno{(*)}
$$
if
\item{1)}
  some vertices and faces are separated out and called {\it
  exterior}, the remaining vertices and faces are called {\it interior};
\item{2)}
    the label of each interior face is a cyclic permutation of one of
  the words $w_i^{\pm1}$;
\item{3)}
  the label of each interior vertex is the identity element of $H$.

A diagram is said to be {\it reduced} if it contains no such
edge $e$ that both faces containing $e$ are interior, these faces
are different and their labels written starting with the $\Mu$-preimages of
$e$ are mutually inverse; such a pair of faces with a common edge is
called a {\it reducible pair}. For example, the faces $C$ and $E$ in Fig.1 
form a reducible pair if  
$\lambda(c_0)=\lambda(e_0)$, $\lambda(c_1)=\lambda(e_2)$, 
$\lambda(c_2)=\lambda(e_1)$ and all edges have the same label.

The following lemma is an analogue of the van Kampen lemma for relative
presentations.

\Lemma 1 {\rm[H83]}.
The natural mapping from a group $H$ to the group with relative
presentation $(*)$ is noninjective if and only if there exists a
spherical diagram over this presentation with no exterior faces and a
single exterior vertex whose label is not 1 in $G$. A
minimal \(with respect to the number of faces\) such diagram is
reduced. \hfil\break
If this natural mapping is injective, then we have
the equivalence: the image of an element $u\in
H*F(t_1,t_2,\dots)\setminus \1$ is 1 in the group $(*)$ if and only
if there exists a spherical diagram over this presentation without
exterior vertices and with a single exterior face with label $u$. A
minimal \(with respect to the number of faces\) such diagram is
also reduced.

Diagrams on the sphere with a single exterior face and no exterior
vertices are also called {\it disk diagrams}, the boundary of the exterior
face of such a diagram is called the {\it contour} of the diagram.

Let $\phi\:P\to P^\phi$ be an isomorphism between two subgroups of a group
$H$.  A relative presentation of the form
$$
\gp{H,t\ |\ \{p^t=p^\phi;\
p\in P\setminus\1\}, w_1=1,\ w_2=1,\ \dots}
\eqno{(**)}
$$
is called a {\it
$\phi$-presentation}. A diagram over a $\phi$-presentation $(**)$ is
called {\it $\phi$-reduced} if it is reduced and different interior cells
with labels of the form $p^tp^{-\phi}$, where $p\in P$, have no common
edges.

\Lemma 2 {\rm[Kl05]}. A minimal \(with respect to the number of faces\) 
diagram among all spherical diagrams over a given $\phi$-presentation 
without exterior faces and with a single exterior vertex with nontrivial 
label is $\phi$-reduced. If no such diagrams exists, then a minimal 
diagram among all disk diagrams with a given label of contour is 
$\phi$-reduced.  \rm In other words, the complete $\phi$-analogue of Lemma 
1 is valid.

The idea of the proof is shown in Fig.2.

\goodbreak
\bigskip
\centerline{\input 2.PIC}
\nobreak%
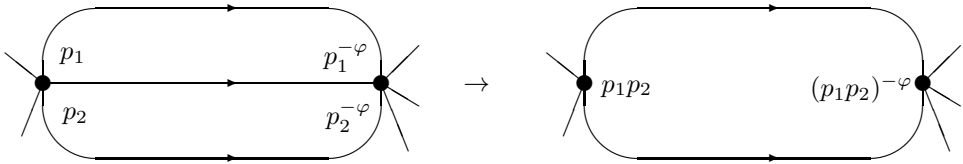
\centerline{Fig. 2}%
\goodbreak
\bigskip

A relative presentation ($\phi$-presentation) over which there exists
no reduced (respectively, $\phi$-reduced) spherical diagrams with 
no exterior faces and
a single
exterior vertex are called {\it aspherical} (respectively, {\it
$\phi$-aspherical}).

\Lemma 3.
Suppose that $H$ is a group, a word
$v\in H*F(t_1,t_2,\dots)$ is not a proper power
in $H*F(t_1,t_2,\dots)$, and a positive integer $l$ is such that
$v^l$
is not conjugate in
$H*F(t_1,t_2,\dots)$ to elements of the set $H\cup \{w_i^{\pm1}\}$
and the presentation
$$
L=\pres<H,t_1,t_2,\dots | v^l=1,w_1=1,w_2=1,\dots>
$$
obtained from presentation $(*)$ by adding the relation $v^l=1$ is
aspherical (or $\phi$-aspherical, if the
initial presentation~$(*)$ is a $\phi$-presentation).
Then
\item{\rm1)}
in the group $K$ with presentation $(*)$, the centraliser of the element
$v^k$ coincides with the cyclic group~$\gp v$ for any positive integer $k$;
\item{\rm2)}
if the group $H$ is nontrivial, then the centre of the group $K$ is
trivial.

\Proof
The first assertion is proven by standard argument. First, we can
assume that $k=l$, because $C(v^k)\subseteq C(v^{kl})$ and the
asphericity of presentation $L$ implies that of the presentation
$$
L_k=\pres<H,t_1,t_2,\dots | v^{kl}=1,w_1=1,w_2=1,\dots>
$$
for each positive integer $k$ (because a cell with label $w^k$ can be
transformed into $k$ cells with labels $w$; see [BoP92]).

Consider a word
$u$ commuting with $v^l$ in the group $K$ and a disk diagram over
presentation $(*)$ with contours
labelled by $[u,v^l]$. Let us glue together segments of the contour of this
diagram to obtain an annulus $A$
with contours labelled by
$v^l$ and $v^{-l}$
(if necessary,
we add cells, whose boundary labels equal 1 in $H*F(t_1,t_2,\dots)$).
Attaching two new cells
$\Gamma_+$ and $\Gamma_-$
to the annulus $A$
along the contours, we
obtain a spherical diagram $D$ over
presentation $L$ without external vertices and faces.
This diagram has the following properties (Fig.3):
\item{a)}
the label of the face $\Gamma_\pm$ written starting with some point
$p_\pm\in\d\Gamma_\pm$ is $v^{\pm l}$;
\item{b)}
the labels of the other faces belong to $\{w_i^{\pm1}\}\cup\1$;
\item{c)}
the points $p_+$ and $p_-$ are joined by a path $\pi$ with label $u$;
\item{d)}
the point $p_+$ is joined with some point $p_-'\in\d\Gamma_-$ by a
path $\pi'$ whose label equals 1 in the free product $H*F(t_1,t_2,\dots)$,
and the label of the
cell $\Gamma_-$ written starting with the point $p_-'$ is $v^{-l}$.

\noindent
The last property follows from the asphericity of the presentation $L$:
in the diagram $D$ and in all diagrams obtained from $D$
by
reductions (eliminations of reducible pairs), the cell $\Gamma_+$ can
form a reducible pair only with the cell $\Gamma_-$.

\goodbreak
\bigskip
\centerline{\input 3.PIC}
\nobreak%
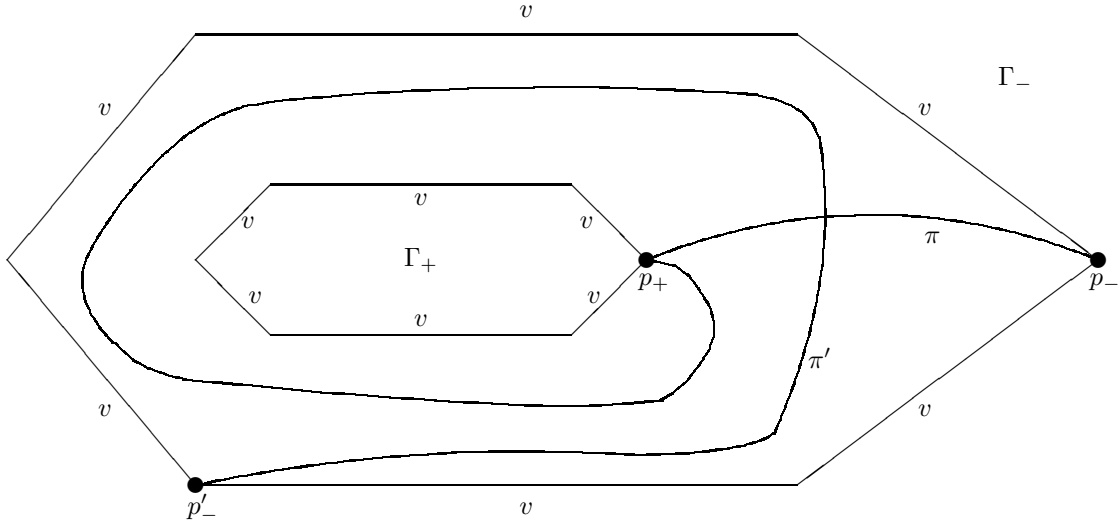
\centerline{Fig. 3}%
\goodbreak
\bigskip

Properties a) and d) and the condition that $v$ is not a proper power
imply that the segment $\sigma$ of the boundary of
$\Gamma_-$ between the points $p_-$ and $p_-'$ has label $v^i$. It
remains to note that the path $\pi\sigma$ is homotopic in the annulus
$A=D\setminus\{\Gamma_\pm\}$ to a path of the form $\pi'\delta^j$, where
$\delta$ is a path with label $v^{-l}$ around the cell $\Gamma_-$
starting and ending at~$p_-'$. This homotopy implies that the label $u$
of the path $\pi$ equals $v^{-lj-i}$ in the group $K$. This proves the
first assertion the lemma.

Let us prove the second assertion. According to assertion 1), the centre
of $K$ must be contained in the cyclic group~$\gp v$. Applying once again
assertion 1) to a hypothetical central element of the form $v^k$, we
obtain $K=\gp v$. Therefore, $v^{ks}=h$ for some $s\in\N$ and
$h\in H$. This contradicts the ($\phi$-)asphericity of the presentation $L$.


\s 2. Proof of Theorem 1

In [Kl05], we showed that the group $\~G$ always (with some obvious
exceptions) admits a relative ($\phi$-)presentation, which is
($\phi$-)aspherical and remains such after adding some additional
relation. By virtue of Lemma 3, this implies Theorem~1. More precisely,
the proof decomposes into two cases.

\noindent
{\bf Case 1: the word $w$ has the form $ct\prod_{i=0}^m(b_ia_i^t)=1$, where
$c,a_i,b_i\in G$} (i.e.,
the complexity of $w$ does not exceeds one in terminology of [FoR05]).

\Lemma 4 {\rm([Kl05], Lemma 23)}.
If $G$ is a torsion-free group and $m\ge0$, then there exists a
$d\in\{2,3\}$ for which the presentation
$$
\gp{G, t\ \Biggm|\ ct\prod_{i=0}^m(b_ia_i^t)=1,\
(a^{t^d}b)^4=1}
$$
is aspherical for any elements $a,b\in G$ such that $a^2\notin\gp{a_m}$
and $b^2\notin\gp{b_0}$.

If $G^2\:=\{g^2\;;\;g\in G\}\not\subseteq\gp{a_m}$ and
$G^2\not\subseteq\gp{b_0}$, then the assertion of Theorem 1 follows from
Lemmata 3 and 4. To complete the proof, it remains to apply the
following simple fact, which we leave to readers as an exercise.

\Lemma 5.
If $G$ is a torsion-free group and $\gp{G^2}$ is a cyclic group, then
the group $G$ is cyclic itself.

\noindent
{\bf Case 2: the word $w$ is not conjugate to a word of the form
$ct\prod_{i=0}^m(b_ia_i^t)=1$} (i.e. the complexity of $w$ is higher than
one in terminology of [FoR05]).

In this case, the assertion of Theorem 1 follows immediately from Lemma 3
and the following two lemmata.

\Lemma 6 {\rm([Kl05], Lemma 2; see also [Kl93], [Fer96])}.
The group $\~G$ has a relative
presentation of the form
$$
\~G\iso\gp{H, t\ \Biggm|\ \{p^t=p^\phi,\ p\in P\setminus\1\},\
ct\prod_{i=0}^m(b_ia_i^t)=1},
\eqno{(1)}
$$
where $a_i,b_i,c\in H$, $P$ and $P^\phi$ are isomorphic subgroups of $H$,
and
$\phi\:P\to P^\phi$ is an isomorphism. The groups $H$, $P$
and $P^\phi$ are free products of finitely many isomorphic copies
of $G$. If the word $w$ is not conjugate in $G*\gp t_\infty$ to a word of
the form $ct\prod_{i=0}^m(b_ia_i^t)$, then the groups $P$ and $P^\phi$ are
nontrivial.

\Lemma 7 {\rm([Kl05], Lemma 10)}. If $G$ is a noncyclic
torsion-free group and $P\ne\1$ in presentation
$(1)$, then there exist elements $a,b\in H$ such that the presentation
$$
\~G\!\!\Bigm/\!\!\nc{a^{t^2}b}\iso
\gp{H, t\ \Biggm|\ \{p^t=p^\phi,\ p\in P\setminus\1\},\
ct\prod_{i=0}^m(b_ia_i^t)=1,\ a^{t^2}b=1 }
$$
obtained from presentation $(1)$ by adding the relator
$a^{t^2}b=1$ is $\phi$-aspherical.

Theorem 1 is proven. In the rest of this paper, we study
generalised unimodular presentations.


\s 3.
Free iterated amalgamated products

This section is an extended version of a similar section of [Kl06b].

Let $\{M_j\,;\,j\in J\}$ be a family of groups. We define a group $M_J$,
a
\emph{[strict] free iterated amalgamated product (FIAP)}
of the groups $M_j$, by induction as follows:
\- if $J=\emptyset$, then we set $M_J=\1$;
\- if $J$ is a finite nonempty set, then a
[strict] FIAP of the family
$\{M_j\,;\,j\in J\}$ is any amalgamated free product of the form
$$
M_J=M_{j_0}\zvezda_{H=H^\phi} M_{J\setminus\{j_0\}},
$$
where $j_0$ is an element of the set $J$, $M_{J\setminus\{j_0\}}$
is a [strict] FIAP of the family $\{M_j\,;\, j\in J\setminus\{j_0\}\}$, and
$\phi:H\to H^\phi$ is an isomorphism between a [proper] subgroup
$H\subseteq M_{j_0}$ and some subgroup
$H^\phi\subseteq M_{J\setminus\{j_0\}}$;
\-if the set $J$ is infinite, then a
[strict] FIAP of the family
$\{M_j\,;\,j\in J\}$ is the direct limit
$$
M_J=\lim\{M_K;\ \hbox{$K$ is a finite subset of $J$}\},
$$
where $M_K$ is a [strict] FIAP of the family groups
$\{M_j\,;\, j\in K\}$; for each pair of finite subsets
$K\subset K'\subset J$, there is a homomorphism $M_K\to M_{K'}$ which is
identity
on the groups $M_j$, where $j\in K$; the direct limit is taken
over this family of homomorphisms.

\Remark.
The definition of strict FIAP requires that the amalgamated subgroup
$H$ is proper only in one factor, in~$M_{j_0}$.
Therefore, any nontrivial group $G$ can be decomposed into a strict
FIAP: $G=\1*G$ (while the trivial group is a strict FIAP
of the empty family of groups). Similarly, if a group $G$ is a
union of a strictly increasing chain of subgroups, i.e.,
$G=\bigcup G_i$, where $G_1\subset G_2\subset\dots$, then $G$
decomposes into a strict FIAP of the groups $G_i$.

Let $I$ be a set, and let $\Omega$ be a
family of subsets of $I$. For each $i\in I$, let $G_i$
be a group, and for each
$\omega\in\Omega$, let $G_\omega$ be a quotient of the free product
$\zvezda\limits_{i\in\omega}G_i$:
$$
G_\omega= \(\zvezda_{i\in\omega}G_i\)\Big/N_\omega.
$$
The natural question arises: under what conditions are the natural mappings
$$
\phi_\omega\:G_\omega\to G_I\:=
\(\zvezda_{i\in I}
G_i\)\Bigg/\!\gp{\!\!\!\gp{\bigcup_{\omega\in\Omega}N_\omega}\!\!\!}
$$
injective? Or under what conditions is the group $G_I$
a free iterated amalgamated product of the groups~$G_\omega$?

The following proposition gives some sufficient condition for this
question to have a positive answer.

\Proposition 1.
Suppose that
$$
N_\omega\cap\(\zvezda\limits_{j\in\omega\setminus\{i\}}G_j\)=\1
\eqno{({**}*)}
$$
for each $\omega\in\Omega$ and each
$i\in\omega\setminus\(\bigcap\Omega\)$. Suppose also that,
for each finite subfamily $F\subseteq\Omega$ with
$|F|\ge2$, there exist elements $\min, \max\in \bigcup F$ such
that
\item{\rm 1)}
the element $\min$ belongs to precisely one set
$\omega_{\min}\in F$;
\item{\rm 2)}
the element $\max$ belongs to precisely one set
$\omega_{\max}\in F$;
\item{\rm 3)}
$\omega_{\min}\ne\omega_{\max}$.
\enditem
Then all of the natural mappings $\phi_\omega\:G_\omega\to G_I$ are
injective
and the group $G_I$ is a free iterated amalgamated product of groups
$G_\omega$. If, in addition,
$$
\(\zvezda\limits_{j\in\omega\setminus\{i\}}G_j\)N_\omega
\not\supseteq G_i
$$
for each $\omega\in\Omega$ and each
$i\in\omega\setminus\(\bigcap\Omega\)$, then the free iterated
amalgamated product is strict.

\Example.
Suppose that $I=\{a,b,c,d,e,f\}$ and
$\Omega=
\left\{
\{a,b,d,e\},
\{b,c,e,f\},
\{d,e,f\}
\right\}$.

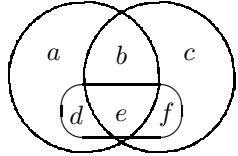
\nobreak
\centerline{\input 0.PIC}
\goodbreak

\noindent
Let $A,\dots,F$ be
the corresponding six groups $G_i$, and let
$\bf ABDE$, $\bf BCEF$, and $\bf DEF$ be the three groups
$G_\omega$.
It is easy to see that conditions 1), 2), and 3) hold for
the family $\Omega$ and each of its two-set subfamilies. Suppose that
condition~$({**}*)$ holds too. Then the validity of Proposition 2 (for
this example) is implied by the following decomposition of $G_I$ into
an amalgamated free product:
$$
G_I= \left( ({\bf DEF}*B)
\zvezda_{B*D*E} {\bf ABDE} \right) \zvezda_{B*E*F} {\bf BCEF}.
$$

\smallskip

To prove Proposition 1 in the general case, we need a lemma.

\Lemma 8 {\rm([Kl06b], Lemma 1)}.
Suppose that the conditions of Proposition 1 hold,
$\Omega'$ is a finite subfamily of $\Omega$,
$\omega\in\Omega$, and
$\alpha\subseteq\omega\cap\left(\bigcup\Omega'\right)$
is a proper subset of $\omega$ contained in
$\bigcup\Omega'$ and containing $\bigcap\Omega$. Then the natural mapping
$$
\zvezda_{i\in\alpha}G_i\to G_{\Omega'}\:=
\left(\zvezda_{i\in\bigcup\Omega'}
G_i\right)\Bigg/
\!\gp{\!\!\!\gp{\bigcup_{\omega'\in\Omega'}N_{\omega'}}\!\!\!}
$$
is injective.

\Proof

\noindent{\bf Case 1: $\omega\in\Omega'$.}
Let us use induction on the cardinality of $\Omega'$.
If $|\Omega'|=1$ (i.e., $\Omega'=\{\omega\}$), then the assertion of
Lemma~8 is true by condition $({**}*)$. Suppose that $|\Omega'|\ge2$. In
this case, according to conditions 1), 2), and 3), the family $F=\Omega'$
contains a set $\omega'\ne\omega$ that has an element $m\in\omega'$
not belonging to $\bigcup(\Omega'\setminus\{\omega'\})$.

By the inductive hypothesis (applied to the set $\omega'$ as $\omega$
and the family $\Omega'\setminus\{\omega'\}$ as~$\Omega'$),
the groups
$$
G_i\ \hbox{with}\
i\in\beta\:=\omega'\cap\left(\bigcup(\Omega'\setminus\{\omega'\})\right)
$$
freely generate their free product in the group
$G_{\Omega'\setminus\{\omega'\}}$.  But according to condition $({**}*)$,
the same groups $G_i$ with $i\in\beta$ freely generate their free product
in the group $G_{\omega'}$ (because $\omega'$ contains an element $m$ not
belonging to $\beta$).  Therefore, the group $G_{\Omega'}$ decomposes into
the amalgamated free product of $G_{\Omega'\setminus\{\omega'\}}$ and
$G_{\omega'}$ with amalgamated subgroup $\zvezda_{i\in\beta}G_i$. The
groups $G_i$ with $i\in\alpha$ lie in the factor
$G_{\Omega'\setminus\{\omega'\}}$. Therefore, the assertion of Lemma 8
follows from the inductive hypothesis applied to the set $\omega$ and the
family $\Omega'\setminus\{\omega'\}$ as $\Omega'$.

\noindent{\bf Case 2: $\omega\notin\Omega'$.}
In this case, the proof is similar.
We again use induction on the cardinality of $\Omega'$.
If $\Omega'=\emptyset$, then we have nothing to prove.
Suppose that $|\Omega'|\ge1$. In this
case, according to conditions 1), 2), and 3), the family
$F=\Omega'\cup\{\omega\}$
contains a set $\omega'\ne\omega$ with an element $m\in\omega'$
not lying in $\bigcup(F\setminus\{\omega'\})$ (see Fig.4).


\goodbreak
\bigskip
\centerline{\input 4.PIC}
\nobreak%
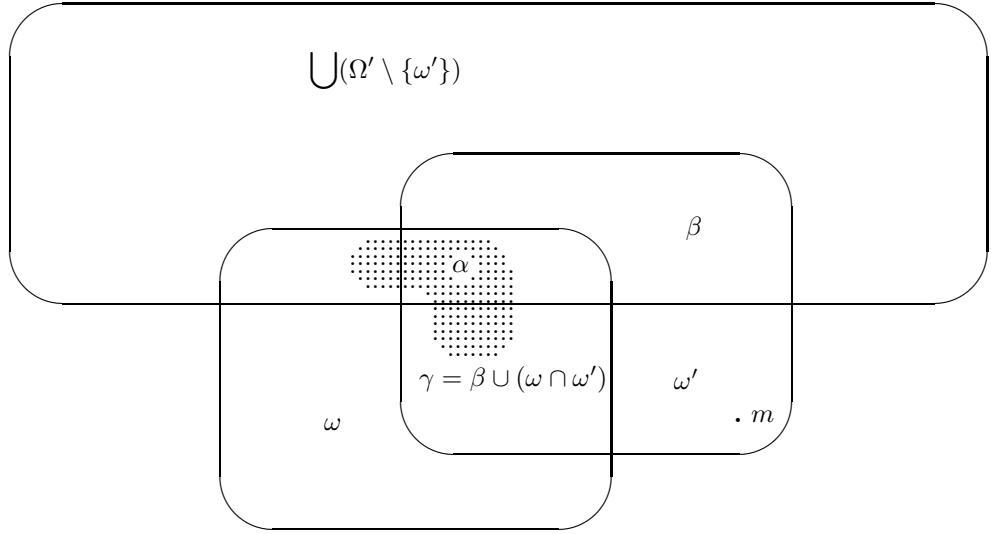
\centerline{Fig. 4}%
\goodbreak
\bigskip

By the inductive hypothesis (applied to the set $\omega'$ as $\omega$ and
the family $\Omega'\setminus\{\omega'\}$ as $\Omega'$), the groups
$$
G_i\ \hbox{with}\
i\in\beta\:=\omega'\cap\left(\bigcup(\Omega'\setminus\{\omega'\})\right)
$$
freely generate their free product in the group
$G_{\Omega'\setminus\{\omega'\}}$.
Therefore, the groups
$$
G_i\ \hbox{with}\
i\in\gamma\:=\beta\cup(\omega\cap\omega')=
\omega'\cap\left(\bigcup((\Omega'\cup\omega)\setminus\{\omega'\})\right)
$$
freely generate their product in the group
$$
H=\left(\zvezda_{j\in(\omega\cap\omega')\setminus\beta}G_j\right)*
G_{\Omega'\setminus\{\omega'\}}.
$$
But condition $({**}*)$ implies that the same groups $G_i$ with
$i\in\gamma$ freely generate their free product in $G_{\omega'}$ (because
$\omega'$ contains an element $m$ not belonging to
$\gamma$). Therefore, the
group $G_{\Omega'}$ decomposes into the amalgamated free product of the
groups $H$ and $G_{\omega'}$:
$$
G_{\Omega'}=H \zvezda_{\gp{G_i\ ;\ i\in\gamma}} G_{\omega'}.
$$
The groups $G_i$ with $i\in\alpha$ lie in
the factor $H$. Therefore, by the inductive hypothesis applied to the set
$\omega$ and the family $\Omega'\setminus\{\omega'\}$ as $\Omega'$, the
groups $G_i$ with
$i\in\alpha\cap\left(\bigcup(\Omega'\setminus\{\omega'\}\right)$
freely generate their free product in
$G_{\Omega'\setminus\{\omega'\}}$. This immediately implies that
the groups
$G_i$ with subscripts $i\in\alpha$ freely generate their free product
in $H$ and, hence, in the group $G_{\Omega'}$, which contains
$H$ as a subgroup. Lemma 8 is proven.

\smallskip\noindent
{\bf Proof of Proposition 1.}
Clearly, it is sufficient to prove Proposition 1 for a finite family
$\Omega$ of cardinality larger than one. In this case,
$$
G_I=G_\Omega*\left(\zvezda_{i\notin\bigcup\Omega}G_i\right),
$$
and $G_\Omega$ decomposes into the amalgamated free product:
$$
G_\Omega=G_{\omega_{\min}}\zvezda_K G_{\Omega\setminus\{\omega_{\min}\}},
$$
where the amalgamated subgroup $K$ is (by virtue of Lemma 8) the free
product of the groups $G_i$ with
$i\in\omega_{\min}\cap\bigcup(\Omega\setminus\{\omega_{\min}\})$.
An obvious inductive argument completes the proof.

In what follows, we need a property of free iterated amalgamated products.

\Proposition 2.
Suppose that a group $M_J$ is a strict free iterated amalgamated product
of finitely generated groups~$M_j$, where $j\in J$.
Then
\item{\rm 1)}
            in the group $M_J$, all the subgroups $M_j$, where $j\in J$,
            are pairwise different;
\item{\rm 2)}
            if an element $h\in M_J$ permutes the groups $M_j$, i.e.,
            for each $j\in J$, there exists a $k\in J$ such that
            $M_j^h=M_k$,
            then $h$ lies in one of the groups $M_i$, where $i\in J$;
            in particular, each central element of $M_J$ is contained 
            in the centre of one of the groups $M_i$.

\Proof
Let us prove the first assertion. Suppose that $M_i=M_k$, where $i$ and $k$
are different elements of the set $J$.
The group $M_i$ is finitely generated; hence, the equality $M_i=M_k$ holds
in a FIAP $M_P$ of a finite family of groups whose direct limit is the
group $M_J$. Here $P$ is a finite set containing $i$ and $k$.
Now, let us use induction on the cardinality of~$P$. The group $M_P$
decomposes into an amalgamated product
$$
M_P=M_{p_0}\zvezda_H M_{P\setminus\{p_0\}}.
$$
The subgroup $H$ is a proper subgroup of $M_{p_0}$ by the definition
of the strict FIAP. Therefore, $M_{p_0}$ can coincide with no
subgroup of $M_{P\setminus\{p_0\}}$. In particular, $i\ne p_0\ne k$.
But then, the equality $M_i=M_k$ holds in the group
$M_{P\setminus\{p_0\}}$. Applying the inductive hypothesis, we complete
the proof of the first assertion.

The second assertion is proven similarly. Suppose that
$h\in\gp{M_{j_1},\dots,M_{j_l}}$. Consider the equalities
$$
M_{j_1}^h=M_{k_1},\dots, M_{j_l}^h=M_{k_l}.
$$
The groups $M_j$ are finitely generated; hence, these equalities hold in
a FIAP $M_P$ of a finite family of groups whose limit is $M_J$.
Here, $P$ is a finite set containing $j_1,\dots,j_l$ and $k_1,\dots,k_l$.
Let us use induction on the cardinality of $P$. The group $M_P$
decomposes into an amalgamated product
$$
M_P=M_{p_0}\zvezda_H M_{P\setminus\{p_0\}}.
$$
The subgroup $H$ is a proper subgroup in $M_{p_0}$ by the definition
of strict FIAP. Moreover, we can assume that $H$ is a proper
subgroup of $M_{P\setminus\{p_0\}}$, because otherwise the element $h$
belongs to $M_{p_0}$ and we have nothing to prove.

If $p_0\notin\{j_1,\dots,j_l,k_1,\dots,k_l\}$, then
$h\in M_{P\setminus\{p_0\}}$, the above equalities hold in the
group~$M_{P\setminus\{p_0\}}$, and the assertion is proven by the
inductive hypothesis.

If $p_0=j_i$, then the group $(M_{p_0})^h$ either lies in
the factor $M_{P\setminus\{p_0\}}$ or coincides with $M_{p_0}$
(depending on whether or not $k_i$ and $j_i$ are equal).
Standard properties of free amalgamated products implies that the
first case is impossible, and the second case is possible only if
$h\in M_{p_0}$.

The case $p_0=k_i$ is considered similarly: we repeat the argument
of the preceding paragraph replacing $h$ by $h^{-1}$.
Proposition 2 is proven.

\Remark.
It can be shown that Proposition 2 becomes false if we remove either of
the condition: the strictness of the FIAP or the finite generatedness of
the factors.


\s 4.
Amalgamated semidirect products

Suppose that a group $A$ acts on a group $B$ by automorphisms
$\phi\:A\to\Aut B$, and $N\nin A$ and $N^\psi\subseteq B$ are
isomorphic subgroups of $A$ and $B$ such that $N$ is normal in $A$
and the isomorphism $\psi\:N\to N^\psi$ is
\emph{consistent with the action $\phi$}, i.e.,
$$
(n^a)^\psi=(n^\psi)^{a^\phi}
\quad\hbox{and}\quad
b^{n^\phi}=b^{n^\psi}
\quad\hbox{for all $a\in A$, $b\in B$, and $n\in N$}.
$$
Consider the semidirect product $A\semitimes B$ corresponding to the
action $\phi$, i.e., the group consisting of the formal products $ab$
multiplied by the formula $(ab)(a_1b_1)=(aa_1)(b^{a_1^\phi}b_1)$.
The elements of the form $nn^{-\psi}$, where $n\in N$, constitute
a subgroup in the semidirect product $A\semitimes B$:
$$
\eqalign{
(nn^{-\psi})(n_1n_1^{-\psi})&=(nn_1)((n^{-\psi})^{n_1^\phi}n_1^{-\psi})=
                              (nn_1)((n^{-\psi})^{n_1^\psi}n_1^{-\psi})=
                              (nn_1)((n^{-1})^{n_1}n_1^{-1})^\psi=
                              (nn_1)(nn_1)^{-\psi};
\cr
(nn^{-\psi})^{-1}&=n^{\psi}n^{-1}=n^{-1}(n^{\psi})^{n^{-\phi}}=
n^{-1}(n^{\psi})^{n^{-\psi}}=n^{-1}(n^{n^{-1}})^\psi=n^{-1}n^\psi.
}
$$
This subgroup is normal:
$$
\eqalign{
(nn^{-\psi})^a&=n^a(n^{-\psi})^{a^\phi}=n^a(n^{-a})^\psi=n^a(n^a)^{-\psi};
\cr
(nn^{-\psi})^b&=b^{-1}(nn^{-\psi})b=nb^{-n^\phi}n^{-\psi}b=
nb^{-n^\psi}n^{-\psi}b=nn^{-\psi}.
}
$$
The {\it semidirect product of the groups $B$ and $A$ (with respect to the
action
$\phi$) with amalgamated (by means of the isomorphism $\psi$)
subgroups $N$ and $N^\psi$} is the quotient group
$$
A\semitimes_{N=N^\psi} B\:=(A\semitimes B)/\{nn^{-\psi}\;;\;n\in N\}.
$$
Clearly, this group contains $A$ and $B$ as subgroups and is the
product of these subgroups; the subgroup $B$ is normal and $A\cap B=N$.


\s 5.
The structure of the group $\^G$

Let $G$ and $T$ be torsion-free groups and let
$w=g_1t_1\dots g_qt_q\in G*T$ be
a cyclically reduced generalised unimodular word, which means that
$$
\prod t_i=t\in\gp{t}_\infty*S=R \nin T
\quad\hbox{and $T/R$ is a group with strong UP property.}
$$
We study the group
$$
\^G=\pres<G,T | w=1>\:=(G*T)\Big/\!\nc{\prod g_it_i}.
$$

First, we describe the main idea of our approach.
The group 
$G*T$ decomposes into the amalgamated semidirect product
$$ 
G*T=T\semitimes_R\nc{G,R}. 
$$
The normal closure $\nc{G,R}$ of $G*R$ in $G*T$ 
decomposes, in its turn, into the free product 
$$
\nc{G,R}=\(\zvezda_{y\in T/R}G^y\)*S*\gp t_\infty,
$$
and 
$$
w\in L_1=\(\zvezda_{y\in X_1}G^y\)*S*\gp t_\infty,
\quad\hbox{where $X_1$ is a finite subset of $T/R$}.
$$
The group $\nc{G,R}$ can be considered as a FIAP of the groups $L_1^x$, 
where $x\in T/R$. The generalised unimodularity guarantees that a  
similar decomposition remains valid for the quotient by the normal
closure of~$w$. The reason is that the quotient of  
$L_1^x$ by the normal closure of $w^x$ has an (ordinary) unimodular 
relative presentation over the group 
$$ 
H_x=\(\zvezda_{y\in xX_1}G^y\)*S.
$$
The details of these arguments are as follows.

Let us decompose $T$ into a union of cosets with respect to the
subgroup $R$:
$$
T=\coprod_{x\in T/R} c_xR, \quad\hbox{where }c_1=1,
$$
and write the relation $\prod t_ig_i=1$ in the form
$$
t\prod_{i=1}^q {g_i}^{c_{x_i}r_i}=1,
\eqno{(2)}
$$
where $r_i\in R$, $x_i=t_it_{i+1}\dots t_qR$,
and $c_{x_i}r_i=t_it_{i+1}\dots t_q$.
Let $X_1$ be the set of all $x_i\in T/R$ that occurs in relation~(2).
For each $x\in T/R$, consider an isomorphic copy $G^{(c_x)}$ of the group
$G$ (assuming that the isomorphism maps each element~$g\in G$ to the
element $g^{(c_x)}\in G^{(c_x)}$). Take also an isomorphic copy
$\=R=\gp{\=t}_\infty*\=S$ of the group $R=\gp{t}_\infty*S$.
We set
$$
H_1=\=S*\(\zvezda_{y\in X_1}G^{(c_y)}\)
$$
and consider the unimodular relative presentation
$$
\~H_1=\pres<H_1,\=t|\=t\prod_i \({g_i}^{(c_{x_i})}\)^{\=r_i}=1>
$$
over the group $H_1$.
The group $\~H_1$ is the quotient of the group
$$
L_1=H_1*\gp{\=t}_\infty=
\=R*\(\zvezda_{y\in X_1}G^{(c_y)}\)
$$
by the normal subgroup $N_1=\nc{\=t\prod{g_i}^{(c_{x_i})\=r_i}}$.
Now, consider the free product
$$
L=\=R*\(\zvezda_{y\in T/R}G^{(c_y)}\)
$$
and the right action $\phi\:T\to\Aut L$ of the group $T$ on $L$:
$$
(\=r)^{x^\phi}           = \={r^x},
\quad
\(g^{(c_y)}\)^{x^\phi} = \(g^{(c_{yx})}\)^{\=a},
\eqno{(3)}
$$
where $x\in T$, $y\in T/R$, and the element $a\in R$ is
uniquely determined by the equality $c_yx=c_{yx}a$.

For each $x\in T/R$, we set
$$
\eqalign{
X_x=X_1x\subseteq T/R,
\quad
&H_x=\=S*\(\zvezda_{y\in X_x}G^{(c_y)}\),
\quad
L_x=L_1^{\chi^\phi}=H_x*\gp{\=t}_\infty=
\=R*\(\zvezda_{y\in X_x}G^{(c_y)}\),
\cr
&N_x=N_1^{\chi^\phi}\nin L_x,
\quad
\~H_x=L_x/N_x,
}
$$
where $\chi\in T$ is any representative of the element $x\in T/R$.

\Proposition 3. The groups $\~H_x$ have the following properties:
\item{\rm 1)}
all of them are isomorphic; an isomorphism $\~H_1\to\~H_x$
acts as follows: $\=t\mapsto \={t^{c_x}}$, $\=s\mapsto \={s^{c_x}}$,
$g^{(c_y)}\mapsto \(g^{(c_{yx})}\)^{\=r}$, where the element $r\in R$
is uniquely determined by the equality $c_yc_x=c_{yx}r$;
\item{\rm 2)}
The group $\~H_x$ has a unimodular relative presentation over
the group isomorphic to $H_1$, which is the
free product of the group $S$ and $p$ isomorphic copies of the group $G$,
where $p=|X_x|=|X_1|$ equals the number of different cosets among
$$
R,\
t_2t_3\dots t_qR,\
t_3t_4\dots t_qR,\
\dots,\
t_qR;
$$
\item{\rm 3)}
in the group $\~H_x$, we have the decomposition
$$
\gp{\=R,\{G^{(c_y)}\ ;\ y\in Y\}}=
\=R*\left(\zvezda_{y\in Y} G_y\right)
\eqno{(4)}
$$
for each proper subset $Y\subset X_x$;
in particular, the natural mapping $\=R\to\~H_x$ is injective
if $w\notin T$;
\item{\rm 4)}
the natural mapping
$$
\~H_x=L_x/N_x\to K\:=L\Bigg/\!\gp{\!\!\!\gp{\bigcup_{y\in T/R}N_y}\!\!\!}
$$
is injective, and the group $K$ is a free iterated amalgamated product
of the groups $\{\~H_x\,;\, x\in T/R\}$;
\item{\rm 5)}
the FIAP $K$ is strict if the group $G$ is noncyclic.

\Proof

\item{\rm 1)}
The isomorphism $\~H_1\to\~H_x$ is the action $\phi$ of the element $c_x$
on $L$. It maps the group~$L_1$ onto the group $L_x$. The subgroup
$N_1\nin L_1$ is mapped onto the subgroup $N_x\nin L_x$, which implies the
isomorphism of the quotients $\~H_1=L_1/N_1$ and $\~H_x=L_x/N_x$.

\item{\rm 2)}
The group $\~H_1$ has the required property by definition. The fulfillment
of this property for the other groups $\~H_x$ follows from the isomorphism
$\~H_1\to\~H_x$ from assertion 1).

\item{\rm 3)}
{%
The group $\~H_1$ has the required property because of the following lemma.
\Lemma 9 {\rm([Kl06a], Lemma 1)}.
Suppose that a group $H$ is torsion-free, a subgroup $P\subseteq H$ is a
free factor of $H$, and a word~$v\in H*\gp z_\infty$ is unimodular and
nonconjugate in $H*\gp z_\infty$ to elements of the subgroup
$P*\gp z_\infty$.
\hfil\break
Then $\nc v\cap (P*\gp z_\infty)=\1$. {\rm In
other words, the element $z$ is {\it transcendental over $P$} in
$\~H=\pres<H,z | v=1>$}.

The other
groups $\~H_x$ have property 3) because of the isomorphism
$\~H_1\to\~H_x$ from assertion 1).
}

\item{\rm 4)}
Let us show that the family of subproducts $\{L_x\ |\ x\in T/R\}$
of the free product $L$, together with the subgroups~$N_x\nin L_x$,
satisfies the conditions of Proposition 1.
Indeed, conditions 1), 2), and 3) of this proposition
follow immediately from the strong UP property of the group
$T/R$.  Condition $({**}*)$ is implied by decomposition (4). Thus,
property~4) is a corollary of Proposition 1.

\item{\rm 5)}
According to Proposition 1, it is sufficient to show that
$G^{(c_y)}\not\subseteq\gp{\=R,\{G^{(c_z)}\,;\,z\in X_x\setminus\{y\}\}}$
for each $x\in T/R$ and each $y\in X_x$ in the group~$\~H_x$.
By virtue of the isomorphism $\~H_x\iso \~H_1$, we can put $x=1$
without loss of generality. So, we have to prove
the impossibility of the inclusion
$$
G^{(c_y)}\subseteq\gp{\=R,\{G^{(c_z)}\,;\,z\in X_1\setminus\{y\}\}}
\quad
\hbox{in the group}
\quad
\~H_1=\pres<H_1,\=t|\=t\prod_i \({g_i}^{(c_{x_i})}\)^{\=r_i}=1>.
\eqno{(5)}
$$
Consider the quotient
$
U=\~H_1/\!\nc{\=S,\{G^{(c_z)}\,;\,z\in X_1\setminus\{y\}\}}.
$
The group $U$ has a unimodular relative presentation over the group
$G^{(c_y)}$. Therefore, the natural mapping $G^{(c_y)}\to U$ is injective
[Kl93]. So, inclusion (5) implies that
the group $G^{(c_y)}$, which is isomorphic to $G$, lies in the cyclic
subgroup $\gp{\=t}$ of $U$. Hence, the group $G$ is cyclic, which
contradicts the assumption. Proposition 3 is proven.

\smallskip\noindent
The action $\phi$ of the group $T$ on $L$ descends to
an action on $K$. We denote this action by the same letter
$\phi\:T\to\Aut K$.
This action is consistent with the isomorphism
$T\supseteq R\iso\=R\subseteq K$.
Indeed, formulae (3) imply that
$$
(\=r)^{x^\phi}=\={(r^x)} \quad\hbox{and}\quad
\(g^{(c_y)}\)^{r^\phi}=\(g^{(c_y)}\)^{\=r}
\quad\hbox{for all $x\in T$, $y\in T/R$ and $r\in R$}.
$$

\Th 3.
The amalgamated semidirect product
$P=T\semitimes\limits_{R=\=R} K$ corresponding to the action $\phi$ and
the isomorphism~$r\mapsto\=r$ is isomorphic to the group $\^G$. The
isomorphism $P\to\^G$ is identity on $T$ and maps the subgroup
$G^{(1)}\subseteq P$ onto $G\subseteq\^G$ and the subgroup $K\nin P$
onto $\nc{G,S}=\nc{G,R}\nin\^G$.

\Proof
The group $G$ embeds in $P$ as a subgroup:
$G=G^{(1)}\subseteq K \subseteq P$. According to the definition of
the action, we have
$
G^{(c_x)}=G^{c_x}.
$
Thus, the relation of the group $\~H_1$, which is valid in $K$, and
the equality $R=\=R$ in the group $P$ give relation (2).
So, the mapping defined by the formulae
$G\ni g\mapsto g^{(1)}\in G^{(1)}$ and $T\ni x\mapsto x$
is a homomorphism from $\^G$ into $P$. The inverse homomorphism has the form
$G^{(c_x)}\ni g^{(c_x)}\mapsto g^{c_x}$ and $T\ni x\mapsto x$.

\Remark. 
Theorem 3 and Proposition 3(3) imply that the natural mapping
$T\to\^G$ is injective if $w\notin T$. In the case when the group $T$ is
free, this fact was first proven by S.~V.~Ivanov, who used
geometrical methods (unpublished).  

\s 6.
Proof of Theorem 2

The first assertion of the theorem follows immediately from Theorem 3.
To prove the second assertion, consider a cyclically reduced
generalised unimodular word $w=g_1t_1\dots g_qt_q$ and put
$t=t_1t_2\dots t_q$. By virtue of Theorem~1, we can assume
without loss of generality that the group $T$ is not cyclic.

\smallskip
\noindent
{\bf Case 1}: $q=1$ and $g_1=1$.
In this case, the centre of the group
$\^G\iso G*(T/\!\nc{t})$
is trivial, if $\nc{t}\ne T$. If  $\nc{t}=T$, then
$R=T$, $S=\1$ (see the definition of unimodularity), and
$T=R=\gp{t}*S=\gp{t}$ which contradicts the assumption that
the group $T$ is noncyclic.

\smallskip
\noindent
{\bf Case 2}: $q=1$ and $g_1\ne1$.
In this case, $T\ne\gp{t}$ (because the group $T$
is assumed to be noncyclic),
$$
\^G\iso
G\zvezda_{\kern 5pt\hbox{$\scriptstyle{g_1=t^{-1}}$}} T,
$$
and the assertion of Theorem 2 is a corollary of the following
well-known simple fact (see, e.g., [LS77]):

\Lemma{10}.
The centre of a free product with amalgamated subgroup which is
proper in both factors coincides with the intersection of the centres
of the factors.

\smallskip
\noindent
{\bf Case 3}:
$q>1$ and the group $\gp{t_1,\dots,t_q}$ is cyclic
(therefore, this group is generated by the element $t=\prod t_i$
by virtue of unimodularity). In this case, the group $\^G$ is
an amalgamated free product:
$$
\^G\iso\gp{G,t\ |\ w=1}\zvezda_{\gp t} T.
$$
By virtue of Lemma 10, the centre of $\^G$ can be
nontrivial only if the centre of $\~G=\gp{G,t\ |\ w=1}$
nontrivially intersects $\gp t$. By Theorem 1, the nontriviality
of the centre of $\~G$ implies that the group $G$ is cyclic, which
contradicts the conditions of Theorem 2.

\smallskip
\noindent
{\bf Case 4}:
the group $\gp{t_1,\dots,t_q}$ is noncyclic. In this case, we can
assume without loss of generality that $T=\gp{t_1,\dots,t_q}$. We
assume also that the group $G$ is noncyclic; we have to prove that
the centre of $\^G$ is trivial.

First, note that it is sufficient to prove the assertion for finitely
generated group $G$. Indeed, the group $\^G$ decomposes into an
amalgamated free product
$$
\^G= G\zvezda_{\gp{g_1,\dots,g_q}}
\pres<\gp{g_1,\dots,g_q}*T|g_1t_1\dots g_qt_q=1>.
$$
If $G\ne\gp{g_1,\dots,g_q}$, then this decomposition implies (by Lemma
10) that the centre of $\^G$ is contained in $G$.
On the other hand, consider the decomposition
$\^G=T\semitimes\limits_{R=\=R} K$ from Theorem 3. According to
the description of the group $K$ (see Section 5), $G$ is contained in
the subgroup $H_1$ of $K$, which is the free product of the group
$S$ and several copies of the group $G$:
$$
G=G^{(1)}\subseteq
H_1=\=S*\zvezda_{y\in X_1}G^{(c_y)}\subset K.
$$
The centre of $H_1$ can be nontrivial only if $S=\1$ and $|X_1|=1$. This
means (by Proposition 3) that $T=R=\gp t$, i.e., the group $T$ is cyclic,
which contradicts the assumption.

In what follows, we assume that $G$ is a finitely generated
noncyclic group and $q\ge2$. We have to prove that the centre of $\^G$
is trivial.
Let us use Theorem 3 again. Suppose that $fy$ is a central element
of the group
$\^G=T\semitimes\limits_{R=\=R}K$, where $y\in T$ and $f\in K$.
According to Proposition 3, the group $K$ is a strict free iterated
amalgamated product of the groups
$\{\~H_x\,;\, x\in T/R\}$. The element $f\in K$ permutes the factors:
$$
\~H_x^f=\~H_x^{y^{-1}}=\~H_{xy^{-1}}.
$$
Therefore, $f\in\~H_z$ for some $z\in T/R$
(by Proposition 2). Thus, the centrality of the element $fy$ means that
$$
\~H_z=\~H_z^{fy}=\~H_z^y=\~H_{zy}.
$$
According to Proposition 2, this equality can hold only for
$z=zy\in T/R$. This implies that $y\in R$ and $fy\in Z(K)$.
Applying Proposition 2 once again, we see that $fy\in Z(\~H_z)$
for some $z\in T/R$. To complete the proof, it remains to
recall that, according to Theorem 1, the centre of $\~H_z\iso\~H_1$ can be
nontrivial only if $H_1$ is cyclic,
which contradicts the assumption about the noncyclicity of the group
$G\iso G^{(1)}\subseteq H_1$.

\REFERENCES

\[B81]
Brodskii S.D.
Anomalous products of locally indicable groups.
Collection: Algebraic systems, 
Ivanov. Gos. Univ., Ivanovo, 1981. P.51--77 

\[B84]
Brodskii S.D.
Equations over groups and one-relator groups
{// Sib. Mat. Zh.} 1984. {T.25}. no.2. P.84--103.

\[BaTa68]
Baumslag G., Taylor T.
The centre of groups with one defining relator
{// Math. Ann.} 1968. {V.175}. P.315--319.

\[BoP92]
Bogley W. A., Pride S. J.
Aspherical relative presentations
{// Proc. Edinburgh Math. Soc. II}. 1992. {V.35}. no.1. P.1--39.

\[CoLy63]
Cohen D. E., Lyndon R. C.
Free bases for normal subgroups of free groups
{// Trans. Amer. Math. Soc.} 1963. {V.108}. P.528--537.

\[CR01]
Cohen M. M., Rourke C.
The surjectivity problem for one-generator, one-relator extensions of
torsion-free groups
{// Geometry \& Topology}. 2001. {V.5}. P.127--142.
See also
arXiv:math.GR/0009101

\[DuH93]
Duncan A.J, Howie J. 
One-relator products with high-powered relator,
in: Geometric group theory 
(G.A.Niblo, M.A.Roller, eds.), 
P.48--74,
Cambridge Univ. Press, Cambridge (1993).

\[FeR96]
Fenn R., Rourke C.
Klyachko's methods and the solution of equations over torsion-free groups
{// L'Enseignment Math\'ematique.} 1996. {T.42}. P.49--74.

\[FoR05]
Forester M., Rourke C.
Diagrams and the second homotopy group
{// Comm. Anal. Geom.} 2005. {V.13}. P.801-820.
See also
arXiv:math.AT/0306088

\[How83]
Howie J.
The solution of length three equations over groups
{// Proc. Edinburgh Math. Soc.} 1983. {V.26}. P.89--96.

\[How87]
Howie J. 
How to generalize one-relator group theory, 
in: Combinatorial group theory and topology 
(S.M. Gersten and J.R. Stallings, eds.), 
53-78, Ann. of Math. Stud., 111, 
Princeton Univ. Press, (1987).

\[How91]
Howie J. 
The quotient of a free product of groups by a single high-powered relator.
III: The word problem
{// Proc. Lond. Math. Soc.} 1991. {V.62}. No.3 P.590-606. 

\[Kl93]
Klyachko Ant. A.
A funny property of a sphere and equations over groups
{// Comm. Algebra}. 1993. {V.21}. P.2555--2575.

\[Kl05]
Klyachko Ant. A.
The Kervaire--Laudenbach conjecture and presentations of simple groups
{// Algebra i Logika}. 2005. {T. 44}. {no.4}. P. 399--437.
See also
{arXiv:math.GR/0409146}

\[Kl06a]
Klyachko Ant.A.
How to generalize known results on equations over groups
{// Mat. Zametki}. 2006. {T.79}. no.3. P.409--419.
See also
arXiv:math.GR/0406382.

\[Kl06b]
Klyachko Ant. A.
The SQ-universality of one-relator relative presentations
{// Mat. Sbornik}. 2006. {T.197}. no.10. P.87--108.
See also
arXiv:math.GR/0603468.

\[Kl07]
Klyachko Ant. A.
Free subgroups of one-relator relative presentations
{// Algebra i Logika}. 2007. {V.46}. no.3. P.290--298
See also 
arXiv:math.GR/0510582.

\[LS77]
Lyndon R. C., Schupp P. E.
{Combinatorial group theory},
Springer-Verlag, Berlin/Heidelberg/New~York, 1977.

\[Met01]
Metaftsis V. 
On the structure of one-relator products of locally indicable groups with 
centre
{// J. Pure Appl. Algebra}. 2001. {V.161}. No.3. P.309--325. 

\[Mu64]
Murasugi K.
The center of a group with a single defining relation
{// Math.Ann.} 1964. {V.155}. P.246--251.

\[Pi74]
Pietrowski A.
The isomorphism problem for one-relator groups with non-trivial centre
{// Math.Z.} 1974. {V.136}. P.95--106.

\[Pr88]
Promyslow S. D.
A simple example of a torsion free nonunique product group
{// Bull. London Math. Soc.} 1988. {V.20}. P.302--304.

\[RS87]
Rips E., Segev Y.
Torsion free groups without unique product property
{// J. Algebra} 1987. {V.108}. P.116--126.

\end